\documentclass{amsart}

\usepackage{latexsym}
\usepackage{oldlfont}
\newtheorem{lemma}{Lemma}[section]
\newtheorem{theorem}[lemma]{Theorem}
\newtheorem{remark}[lemma]{Remark}

\newcommand{\EX}{{\Bbb{E}}}
\newcommand{\PX}{{\Bbb{P}}}

\newcommand{\om}{\omega}

\renewcommand{\phi}{\varphi}

\newcommand{\R}{{\mathbb R}}

\input amssym.def
\input amssym
\title
[Stability of the Stochastic Quasigeostrophic Equation]
{Exponential Stability of the Quasigeostrophic Equation under Random
Perturbations}
\author{Jinqiao Duan}
\address[Jinqiao Duan]
{Department of Applied Mathematics\\
Illinois Institute of Technology\\
Chicago, IL 60616, USA}
\email
[Jinqiao Duan]{duan@math.clemson.edu}

\author{Peter E. Kloeden}
\address[Peter  E. Kloeden]
{Department of Mathematics\\
Johann Wolfgang Goethe University\\
D--60054  Frankfurt am Main,
Germany}
\email [
P.~E.~Kloeden]{kloeden@math.uni-frankfurt.de}

\author{Bj{\"o}rn Schmalfu{ss}}
\address[Bj{\"o}rn Schmalfu{ss}]
{Department of Applied Sciences\\
University of Technology and Applied Sciences\\
Geusaer Stra{ss}e\\
D--06217 Merseburg\\
Germany\\}
\email
[Bj{\"o}rn Schmalfu{ss}]{schmalfuss@in.fh-merseburg.de}

\date{January 1, 1999}


\begin{document}

\maketitle

\begin{abstract}

{\em Progress in Probability} {\bf 49}(2001), 241-256.

The quasigeostrophic  model describes
large scale and relatively slow fluid motion in geophysical flows.
We investigate the quasigeostrophic  model under random forcing and
random boundary conditions. We first transform the model into
a partial differential equation with random coefficients.
Then we show that, under suitable  conditions on the random
forcing, random boundary conditions, viscosity, Ekman constant
and Coriolis parameter,   all quasigeostrophic motion
approach a unique stationary
state exponentially fast.
This stationary state corresponds to a
unique invariant Dirac measure.

\end{abstract}

\section{Introduction}

The quasigeostrophic (QG) model is a simplified geophysical fluid model at
asymptotically high rotation rate or at small Rossby number.  It is derived as an
approximation of the rotating shallow water equations by a conventional
asymptotic expansion for small Rossby number \cite{Ped87}.  The lowest order
approximation gives the barotropic QG equation, which is also the conservation
law for the zero--th order potential vorticity.  Warn et al.  \cite{WaBoShVa95}
and Vallis \cite{Val96} emphasize that this asymptotic expansion is generally
secular for all but the simplest flows and propose a modified asymptotic method,
which involves expanding only the fast modes.  The barotropic QG equation also
emerges at the lowest order in this modified expansion.

Moreover, it has
recently been shown \cite{Schoc87, BeaBou94, EmbMaj96} that quasigeostrophy is a
valid approximation of the rotating shallow water equations in the limit of zero
Rossby number, i.e., for asymptotically high rotation rate.  The
three-dimensional baroclinic quasigeostrophic flow model can be derived
similarly; see, for example, \cite{Ped87, EmbMaj96, BeaBou94, Holm, Desjardins1}.

We consider the barotropic quasigeostrophic
flow model \cite{Ped87, Ped96, Mul96}
\begin{equation}
   \Delta \psi_t + J(\psi, \Delta \psi ) + \beta \psi_x
   =\nu \Delta^2 \psi - r \Delta \psi + {\dot{W}_2}   \label{qg}
\end{equation}
on  a     rectangle $D=(0,1)\times(0,1)$ $\subset$ $\R^2$ ,
where $\psi(x,y,t)$ is the stream function,  $\beta$ $\geq$ $0$  the
meridional
gradient of the Coriolis parameter, $\nu$ $>$ $0$  the viscous dissipation
constant,
$r$ $>$ $0$  the Ekman dissipation constant, ${\dot{W}_2}$  the
 noise due to wind
forcing,
and $J$  the Jacobian operator, which is defined by $J(f, g)$ $=$ $f_x g_y
-f_y g_x$.
Equation (\ref{qg}) can be rewritten in terms of  the vorticity  $q$ $=$
$\Delta \psi$  as
\begin{equation}
       q_t + J(\psi, q ) + \beta \psi_x
        =\nu \Delta q - r q  + {\dot{W}_2} \; ,   \label{eqn}
\end{equation}
which is usually supplemented with boundary conditions
\begin{eqnarray}
 \psi(x,y,t) & = & 0 \quad \mbox{on} \; \partial D \; , \label{BC1} \\
\frac{\partial}{\partial n} q (x,y,t) & = & \dot{W_1} \quad \mbox{on} \; \partial D \; , \label{BC2}
\end{eqnarray}
along with  an appropriate initial condition
\begin{eqnarray}
q (x,y,0)  =  q_0(x,y)                    \;  ,
\end{eqnarray}
where $n$ is the unit outward normal vector on the
boundary $\partial D$,  $W_1$ is a $temporally$ {\em two--sided} Wiener process
with values in the function space $U$.
The boundary condition (\ref{BC1}) means no normal flow can pass through
the boundary. The other boundary condition (\ref{BC2}) says that
$q=\Delta \psi$ has zero mean but with fluctuations, and thus
 might  be called  a
random slip boundary condition \cite{Ped96}. As discussed in   Pedlosky's book
(\cite{Ped96},
page 34), the boundary condition $\Delta \psi =0$ may be
an appropriate slip boundary condition for the
large scale quasigeostrophic motion. Boundary conditions for
the quasigeostrophic model are not quite well
understood, since this model describes large scale flows while boundary
conditions
also involve small scale motions. For this reason we believe that, under
the random
media or random wind forcing conditions,  a random slip condition may be
more
appropriate than the usual slip boundary
conditions \cite{Ped96} for the deterministic quasigeostrophic model.
We also note that the Neumann form for the boundary condition
on $q$ is for mathematical convenience.

In this article  we treat the  quasigeostrophic flow model with both random
forcing and random boundary condition  as a mathematical random dynamical
system \cite{Arn98}.  Our aim is to show that there exists a random steady
state
under a particular choice of parameter values. This random steady state  is a
statistically stationary solution  towards which  any other solution
trajectory
tends as $t\to\infty$. Our mathematical approach  is to formulate the
random quasigeostrophic flow model as   a stochastic evolution equation  with
structural similarities to  the Navier--Stokes equation and then to show
that it
generates a random dynamical system for which there  exists an attracting
 random fixed point.

\section{Preliminaries}

 Following Arnold \cite{Arn98}  we  will model  noise in an abstract  random
dynamical system on a state space $H$ by a   metric dynamical  system
$(\theta,{\cal P})$ on a  probability space   ${\cal P}$ $=$ $(\Omega,{\cal
F},\PX)$.
A {\em metric dynamical system} consists of a group
$\{\theta_t\}_{t\in\Bbb{R}}$
of  operators $\theta_t$ $:$ $\Omega$  $\mapsto$ $\Omega$,  i.e.,  satisfying
\[
\theta_0 = id_{\Omega}, \qquad \theta_{s+t}=\theta_s\circ\theta_t \quad
\mbox{for all}
\, \, s,t \in \R,
\]
such that the mapping  $(t,\omega)$ $\mapsto$ $\theta_t\omega$  from $\R
\otimes
\Omega$  into $\Omega$ is $({\cal B}(\Bbb{R})\times{\cal F},{\cal
F})$--measurable
and the probability measure $\PX$ is ergodic (hence invariant) with respect
to the
flow $\theta$.\\

A {\em random dynamical system} consists of  a metric dynamical
system $(\theta,{\cal P})$ and a  cocycle mapping $\phi$ $:$ $\R \times
\Omega \times H$
$\mapsto$ $H$,   i.e.,  satisfying
\[
\phi(0,\omega,\cdot)={\rm id}_H, \quad
\phi(s+ t,\omega,\cdot)=\phi(s,\theta_t,\cdot)\circ\phi(t,\omega,\cdot)
\quad \mbox{for all} \, \, s,t \in \R^+,
\]
that is $({\cal B}(\Bbb{R})\otimes{\cal F}\otimes{\cal B}(H),
{\cal B}(H))$-measurable. The mapping $\phi$  describes the dynamics of the
system in the state space $H$, which will be a separable Hilbert space
with inner product $(\cdot,\cdot)$ and norm  $|\cdot|$ $=$
$\sqrt{(\cdot,\cdot)}$
 in this article. \\

Let $X$ be a random variable defined on $(\theta,{\cal P})$ with values in
$H$. By the invariance of $\PX$ and the measurability of $\theta$ the
mapping $(t,\omega)$ $\mapsto$ $ X(\theta_t\omega)$ from $\R \times \Omega$
 into $H$ is a measurable stationary stochastic process. We will restrict
attention
here to random variables  generating stationary processes that satisfy
certain
growth conditions.  A $H$--valued random variable $X$ is said to be
{\em tempered} with respect to a $\theta$-invariant set $\Omega^\prime$
(of full $\PX$-measure) if the mapping  $t$ $\mapsto$ $ |X(\theta_t\omega)|$
grows at most subexponentially as $t$ $\to$ $\pm\infty$, i.e.,  for which
\[
\lim_{t\to\pm\infty}\frac{\log^+|X(\theta_t\omega)|}{|t|}=0.
\]
for $\omega$ $\in$ $\Omega^\prime$.  Note that  the only alternative to
this  when
 $X$ is not tempered is
\[
\limsup_{t\to\pm\infty}\frac{\log^+|X(\theta_t\omega)|}{|t|}=+\infty.
\]

A {\em random fixed point}  of a random dynamical system $\phi$
is a $H$--valued random variable $X^\ast$  for which
\begin{equation}\label{eq2.1}
\phi(t,\omega,X^\ast(\omega))=X^\ast(\theta_t\omega) \qquad \mbox{for all}
\, \, t \in \R^+
\end{equation}
and  for all $\om$ in  a $\theta$-invariant set of full $\PX$--measure.
The dynamics thus follows a stationary regime if we start in $X^\ast(\omega)$.
In particular, the probability distribution of these states is  independent
of $t$.
\\
The following  theorem is a special case of a random
fixed point theoremdue toSchmalfu{ss} \cite{Schm97a}.

\begin{theorem}\label{theo2.1}
Let $\phi$ be a random dynamical system with a separable Banach space
$(H,|\cdot|)$ as its state space.
Assume that     the mapping $x$ $\mapsto$$\phi(t,\omega,x)$
is continuous for every $t$ $\ge$ $0$ and $\omega$ $\in$  $\Omega$.
In addition,  let   ${\cal X}$ $=$ $\{{\cal
X}(\omega)\}_{\omega\in\Omega}$
be  a closed random set  in $H$ such that
\begin{equation}\label{eq2.2}
\EX \sup_{h_1\not= h_2 \in {\cal X}(\omega) } \log\frac{|\phi(1,\omega,h_1)
-\phi(1,\omega,h_2)|}{|h_1- h_2|} < 0
\end{equation}
and   let the real valued random variable defined by $\omega$ $\to$
$\sup_{h\in {\cal X}(\omega)}|h|$  be tempered.   Moreover, assume that

\begin{equation}\label{eq2.2.1}
\sup_{t\in[0,1]}\sup_{h_1\not= h_2 \in {\cal X}(\omega) } \frac{|\phi(t,\omega,h_1)
-\phi(t,\omega,h_2)|}{|h_1- h_2|}
\end{equation}
is tempered with respect to $\{\theta_n\}_{n\in\mathbb{Z}}$.
\\
Then there exists a $\theta$-invariant set of full $\PX$--measure
$\Omega^\prime$
and a random variable $X^\ast$ satisfying (\ref{eq2.1})  on $\Omega^\prime$.
Futhermore,  $X^\ast$   is exponentially attracting, i.e.,
\[
\lim_{t\to\infty}|\phi(t,\omega,X(\omega))-X^\ast(\theta_t\omega)|=0,\qquad\text{a.s.}
\]
exponentially fast for any measurable selection $X$ of ${\cal X}$.
\end{theorem}

\subsection{Linear stochastic evolution equations}

  In the following we will consider the motion relative to a spatially
constant flow. By the particular structure of the coefficients
this spatially constant flow can by calculated separately. To find
this spatially constant flow we have to solve a simpler equation.
Henceforth we take  for $H$  the space of square integrable functions $L_2(D)$
on the   rectangle $D=(0,1)\times (0,1)$  in $\Bbb{R}^2$   fulfilling
$\int_Df(x)dD=0$
and  denote its  norm by $|\cdot|$.  We then define $V$ with norm
$\|\cdot\|$ to be the Sobolev space of functions   contained in
$H$  with generalized derivatives
of first order belonging to  $L_2(D)$,
and  define $W_2^2(D)$ to  be the space of functions with first and second
generalized derivatives belonging to  $L_2(D)$.  In addition, we denote by
$U=L_2(\partial D)$  a boundary  space associated with   square
integrable functions on the boundary of $D$.
To use the results of DaPrato and Zabczyk \cite{DaPZab96} Chapter 13
we use that this perturbation is only defined on one side of the rectangle,
say $\{0\}\times (0,1)$. So $U$ consists of $L_2$-functions on $\partial D$
which are zero outside of $\{0\}\times (0,1)$ with zero average. However,
generalizations are possible. \\

Let $\Delta$  be the Laplacian operator on  $D$.  The boundary value problem
\[
-\nu\Delta u = f,  \qquad \frac{\partial}{\partial n}u=g\quad\text{on }\partial D
\]
with $f$ $\in$ $H$, $g$ $\in$  $L_2(U)$ and $\nu>0$,
has a unique solution $u$ $=$ $\tilde G(f,g)$. The  solution operator
$\tilde G$ $:$
$H\times  U$ $\to$ $W_2^2(D)$ is a bounded linear
operator,   i.e.,
there exists a constant  $c_G$ such that
\[
\|\tilde G(f,g)\|_{W_2^2(D)} \le c_G \left(|f|+\|g\|_{U}\right).
\]
     Similarly we can consider the same equation but
with the homogeneous Dirichlet boundary condition
\begin{equation}\label{bvp2}
-\Delta u = f,  \qquad u|_{\partial D}=0.
\end{equation}
The solution operator $G$ for this boundary value problem satisfies
\[
\|G(f)\|_{W_2^2(D)}\le c_G |f|.
\]
In the following we will denote by $C_{G,x}$the constant which estimates
$G(f)_x$ with respect to $L_2(D)$-norm of $f$.

In order to introduce  a {\em white noise} on the boundary $\partial D$,
we   consider
 a $temporally$ {\em two--sided} Wiener process $W_1$  with values in
$U$
and denote by  $\{{\mathcal{F}}_t\}_{t \in \R}$  the filtration
corresponding to this  Wiener process,
where,  roughly speaking,  ${\mathcal{F}}_t$   is generated by the
increments of the noise
sample paths between times
$-\infty$ and $t$.  The associated Wiener measure  $P$ is  defined on the
$\sigma$-algebra of the
canonical  sample  space $\Omega$  that  consists   of continuous functions
$\omega$  from $\Bbb{R}$
into a phase space of the noise satisfying $\omega(0)$ $=$ $0$.
We will assume that the
covariance
operator $Q^{W_1}$ with respect to this measure satisfies
${\rm tr}_{U}Q^{W_1}$  $<$ $\infty$   and define
the Wiener shift by
\[
{W_1}(\theta_t\omega,\cdot)={W_1}(\omega,t+\cdot)-{W_1}(\omega,t).
\]
The measure $P$ is ergodic with respect to the flow $\theta$  of the metric
dynamical system
formed by the Wiener shift.

Since the solution operator $\tilde G(0,\cdot)$ above  is a linear and bounded
operator,  the
process $t$ $\mapsto$ $\widetilde {W_1}(\omega,t)$ $:=$ $\tilde G(0,{W_1}(\omega,t))$ also
defines a Wiener process
with trajectories in the Sobolev space $W_2^2(D)$  for which the
bounded covariance operator $\widetilde Q$ $:=$ $\tilde G\circ Q^{W_1}\circ \tilde G^\ast$
has a finite trace with
respect to $W_2^2(D)$. \\

The operator $-\nu\Delta$ with the vanishing Neumann  boundary condition can
be extended to an
operator $A$ defined on $D(A)$ $=$ $W_2^2(D)$   with the
vanishing Neumann boundary condition.
  The space $H$  has
a complete
orthonormal base,   consisting  of eigenvectors $e_1$, $e_2$, $\cdots$ with
corresponding
eigenvalues $\lambda_1$ $\le$ $\lambda_2$ $\le$ $\cdots$ for
the   operator $A$.
By the particular choice of $H$ we know that $\lambda_1>0$ and
 that $A$ is coercive.
  Let the    semigroup $\{S(t)\}_{t\ge 0}$ on $H$   be
the solution operator (indexed by $t$) of
the initial-boundary value problem
\[
\frac{du}{dt}-\nu\Delta u = 0, \qquad  u(0) = u_0 \in H, \qquad
\frac{\partial u}{\partial n} = 0\quad \text{on }\partial D .
\]
with   $\nu$ $>$ $0$.
  This semigroup has the  generator $-A$.

We consider an expression of the form
\begin{equation}\label{eq1.2}
z(\omega,t)=S(t)z_0+ \int_0^tAS(t-\tau)d\widetilde {W_1}(\omega,\tau),
\end{equation}
as the solution of the linear stochastic evolution equation
\begin{equation}\label{eq1.1}
\frac{dz}{dt}-\nu\Delta\,z=0,\qquad z(0)=z_0\in H,\qquad \frac{\partial}{\partial n}z(t)|_{\partial
D}={\dot{W}_1}(t) ;
\end{equation}
see (\cite{DaPZab96} Section 13.2).  The expression (\ref{eq1.2}) is
meaningful if,  for example,
\[
\int_0^t\|AS(t-\tau)GQ^{{W_1}\,\frac{1}{2}}\|_{{\mathcal{L}}^2(U,H)}^2
d\tau<\infty,
 \quad \mbox{for all} \, \, 0 \le t  <+\infty.
\]
Since by  the invariance of the increments of the Wiener process, we then
have
\[
\EX\left|A\int_{-\infty}^0S(-\tau)d\widetilde {W_1}\right|^2
\le
\sum_{i=0}^\infty e^{-2\lambda_1i}\EX\left|A\int_{-1}^0S(\tau)d\widetilde
{W_1}\right|^2<\infty.
\]
The random
variable
$z_{W_1}$ defined  by
\[
z_{W_1}(\omega):=A\int_{-\infty}^0S(-\tau)d\widetilde {W_1}(\omega,\tau)
\]
 is thus well defined and has  finite second moment with respect to the
norm of $H$.  Moreover, we have formally
\begin{eqnarray*}
S(t)z_{W_1}(\omega)&+&A\int_0^tS(t-\tau)d\widetilde {W_1}(\omega,\tau)\\
&=&S(t)A\int_{-\infty}^0S(-\tau)d\widetilde {W_1}(\omega,\tau) +
A\int_0^tS(t-\tau)d\widetilde {W_1}(\omega,\tau)\\
&=& A\int_{-\infty}^0S(t-\tau)d\widetilde
{W_1}(\omega,\tau)+A\int_0^tS(t-\tau)d\widetilde
{W_1}(\omega,\tau)\\
&=&A\int_{-\infty}^tS(t-\tau)d\widetilde {W_1}(\omega,\tau)\\
&=&  A\int_{-\infty}^0S(-\tau)d\widetilde
{W_1}(\theta_t\omega,\tau)=z_{W_1}(\theta_t\omega),
\end{eqnarray*}
so  the stationary process
\[
t \mapsto z(\theta_t\omega) = A\int_{-\infty}^0S(-\tau)d\widetilde
{W_1}(\theta_t\omega,\tau),
\quad  t\in\Bbb{R}
\]
solves the boundary value problem (\ref{eq1.1}).

Since $\EX|z_{W_1}|^2$ is finite,  we can apply the Burkholder inequality to
obtain
\[
\EX\sup_{t\in [0,1]}|z_{W_1}(\theta_t\omega)|^2<\infty,
\]
and it   then follows  from  the Birkhoff Ergodic Theorem \cite{DaPZab96} that
\[
\lim_{i\to\pm\infty}\frac{\sup_{\tau\in
[0,1]}|z_{W_1}(\theta_{i+\tau}\omega)|^2}{i}=0
\]
on a $\theta$-invariant subset of $\Omega$ of full $P$--measure.
Hence
\[
\lim_{t\to\pm\infty}\frac{|z_{W_1}(\theta_t\omega)|^2}{t}=0
\]
on a $\theta$-invariant subset of $\Omega$ of full $P$--measure, i.e.,
$|z_{W_1}|$ is tempered.
Note that  similar techniques  can be used to show  that $z_{W_1}$ is defined on a
$\theta$--invariant
set of full measure.\\

Finally,   equations for the generalized spatial derivatives of  $z_{W_1}$ can
be investigated   if we suppose that the covariance $Q^{{W_1}}$
is sufficiently regular. Conditions are formulated in
DaPrato and Zabczyk \cite{DaPZab96} Theorem 13.3.1.
In particular, $\nabla z_{W_1}$  is well defined   and tempered.

\section{Transformation of the quasigeostrophic equation}

We return to the QG vorticity equation (\ref{eqn}),  in which we now write
$\Delta\psi$ for
the vorticity.  That is, we consider
\begin{equation}\label{eq3.2}
\frac{d\Delta\psi}{dt}+J(\psi,\Delta\psi)+\beta\psi_x=\nu\Delta^2\psi-r\Delta\psi+{\dot{W}_2}(\omega,t),
\end{equation}
with a non zero boundary condition
\begin{equation}\label{eq3.4}
\frac{\partial}{\partial n}\Delta\psi(t)={\dot{W}_1}(\omega,t), \quad\psi=0\quad \text{on }\partial D,
\end{equation}
that involves  a white noise  ${\dot{W}_1}$  on the boundary $\partial D$ as
described in the
previous section. In addition,  the wind forcing white noise $\dot{W}_2$ is
based on
 a temporally two--sided  noise adapted Wiener process $W_2$ with values in
$V$ and
covariance $Q^{W_2}$ such that ${\rm tr}_{V} Q^{W_2}<\infty$.
In particular,  $W_1$ and $W_2$ are assumed to be
independent.

We can  now  define a metric dynamical  system with  the properties of our
white noise
terms.  For $\Omega$ we  choose an appropriate subset of  the function space
$C_0(\Bbb{R},U)$ $\times$ $ C_0(\Bbb{R},V)$
with the usual Borel
$\sigma$--algebra of a Fr\'echet space, i.e., an element  $\omega$ is  a
continuous path
from  $\Bbb{R}$  into  $(U,V)$ with $\omega(0)=0$.
Then we take $\mathbb{P}$ $=$ $P^{W_1}\otimes P^{W_2}$  to be the product measure of the
Wiener measures
 corresponding to  $W_1$ and $W_2$,  which is ergodic  since both $P^{W_1}$ and
$P^{W_2}$ are ergodic.
The flow $\theta$ on $\Omega$ is defined in terms of  shift operators
$\theta$ applied
to the sample paths of $W_1$ and $W_2$. \\

The above QG equation has
structural
similarities to  equations of Navier--Stokes type.  To  be able to  adapt
well  known
results of such equations,   we need to replace these boundary conditions
by zero
boundary conditions, which is possible with particular types of  stationary
transformations;
see Crauel and Flandoli \cite{CraFla94} or Brannan, Duan and Wanner
\cite{BraJinWan98},
 or in a more general context Keller and Schmalfu{ss} \cite{KelSchm97b} or
Imkeller and
Schmalfu{ss} \cite{ImkSchm98}.  In particular, we  transform (\ref{eq3.2}) into
\begin{eqnarray}  \label{eq3.3}
&&
\frac{d}{dt}u + J(G(u),u)+\beta G(u)_x   =
\nu\Delta u-r\,u+{\dot{W}_2}(\omega,t)\\
&&
\qquad  \frac{\partial}{\partial n}u(t)|  =  {\dot{W}_1}(\omega,t)\quad\text{on }\partial D\nonumber
\end{eqnarray}
where  $G$ is the solution operator of the boundary value problem (cf.
(\ref{bvp2}))
\[
\Delta\psi = u, \qquad \psi|_{\partial D} = 0,
\]
i.e., with  the  solution  $\psi$ $=$ $G(u)$.
  We consider equation (\ref{eq3.3}) as an evolution equation on the triple
$V\subset H\subset V^\prime$ where $V^\prime$ is the dual space of $V$.
\\

The properties of the nonlinear term  of  equation (\ref{eq3.3})
follow from those of the  bilinear operator $B$ $:$ $L_2(D)\times W_2^1(D)$
$\to$
$V^\prime$ defined by
\begin{equation}\label{eq3.1}
B(v_1,v_2) = J(G(v_1),v_2).
\end{equation}

\begin{lemma}\label{lem2.3}  \quad $B$  is a well defined,   continuous
operator
and
\begin{eqnarray*}
i)&& \qquad \langle B(v_1,v_2),v_3\rangle=-\langle B(v_1,v_3),v_2\rangle,
\\
ii)&&\qquad \langle B(v_1,v_2),v_2\rangle=0.
\end{eqnarray*}
for $v_1\in L_2(D)$, $v_2\in W_2^1(D),\,v_3\in  V$
\end{lemma}
\begin{proof}   There exist  positive constants $c$, $\,c^\prime$ and
$\,c_B$  such that for
any  $v_1$ $\in$ $L_2(D)$, $v_2$ $\in$ $W_2^1(D)$ $\subset$ $L_4(D)$ and
$v_3$ $\in$
$V$ $\subset$ $L_4(D)$ we have
\begin{eqnarray*}
\qquad|\langle B(v_1,v_2),v_3\rangle|&=&
\left|\int_D(G(v_1)_xv_{2y}v_3-G(v_1)_yv_{2x}v_3)dD\right|\\
&\le&
c\|\nabla G(v_1)\|_{L_4(D)}|\nabla v_2|\|v_3\|_{L_4(D)}\\
&\le& c^\prime \|\nabla G(v_1)\|_{W_2^1(D)}|\nabla v_2|\|v_3\|\\
&\le& c_B |v_1||\nabla v_2|\|v_3\|,
\end{eqnarray*}
 which implies that $B$ is  well defined and continuous.\\ Property
$i)$  follows from the  integration by parts formula:
\begin{eqnarray*}
&&
\int_D G(v_1)_xv_{2\,y}v_3dD-\int_D G(v_1)_yv_{2\,x}v_3dD\\
&&
\quad=-\int_D G(v_1)_xv_{3\,y}v_2dD+\int_D G(v_1)_yv_{3\,x}v_2dD\\
&&
\qquad
+\int_{\partial D}G(v_1)_xv_2v_3\cos(n,y)dS-
\int_{\partial D}G(v_1)_yv_2v_3\cos(n,x)dS\\
&&\quad=-\langle B(v_1,v_3),v_2\rangle
\end{eqnarray*}
  because the boundary integrals are zero. Indeed, for two
sides of $\partial D$ these integrals are zero by the orthogonality
of $n$ and the direction of the derivative. For the other both sides
the integrals are also zero. For example, for the first
integral we have by the properties of $G$
\[
G(v_1)[x,1]=G(v_1)[x,0]=0,
\]
hence $G(v_1)_x[x,1]=G(v_1)_x[x,0]=0$ .\\
Property
$ii)$ is a consequence of the  antisymmetric nature of  property
$i)$.
\end{proof}

\begin{remark}
{\rm    If $v_1\in V$,  which one can assume to be the solution of
(\ref{eq3.3}), then we can similarly  get
 that $B(v_1,v_1)\in H\subset L_2(D)$. This shows that we can
split up the solution of the original equation into a
special constant part plus the remaining part.
Similarly we get $G(v_1)_x\in H$.}
\end{remark}

Equation (\ref{eq3.3}) is similar to  the equations of the Navier--Stokes type. Indeed, the Laplace operator term  in
(\ref{eq3.3})
is also present in the Navier--Stokes equations (see Temam \cite{Tem79}),
while
 the bilinear operator $B$ defined by  (\ref{eq3.1}) has similar properties
(actually, a bit stronger) to  the bilinear operator defining
the nonlinearity of   the $2$--dimensional  Navier--Stokes equations.
It  thus  follows from the general theory of the stochastic Navier--Stokes
equation
that (\ref{eq3.3}) has a unique solution, see for instance Schmalfu{ss}
\cite{Schm95d}.
The linear terms $ru$ and $\beta G(u)_x$ appearing in (\ref{eq3.3}) but not
in the Navier--Stokes equation are not essential for a proof of existence and
uniqueness.  See Brannan, Duan and Wanner \cite{BraJinWan98}  for another
proof of
existence and uniqueness based on mild solutions.

\section{The stationary solution}

We now transform the stochastic evolution equation (\ref{eq3.3})
into a random
evolution equation in $V$ $\subset$ $H$ $\subset$ $V^\prime$,  i.e.,
with
 stationary random coefficients rather than white noise driving or boundary
terms.
This  will make it easier to find a  forward invariant random set on which
we can verify
an appropriate Lipschitz  condition.  We introduce  the random variable
\[
z_{W_2}(\omega)=\int_{-\infty}^0S(-\tau)dW_2(\omega,\tau)\in V,
\]
which we note without proof is  a tempered random variable on a
$\theta$--invariant set of full measure.
We also assume that $W_2$ (hence $Q^{W_2}$) is sufficiently regular
such that the Neumann boundary condition is fulfilled.
  Since $z_{W_2}$ fulfills the Neumann boundary condition  there is no influence
to the boundary condition  of (\ref{eq3.3}). \\
We consider the  random evolution equation
\begin{equation}\label{eq4.0}
\begin{split}
\frac{d}{dt}z(t,\omega)&+B(z,z)+ Az+\beta\,G(z)_x+rz
\\
=&-B(z,z_{W_1}(\theta_t\omega)+z_{W_2}(\theta_t\omega))-
B(z_{W_1}(\theta_t\omega)+z_{W_2}(\theta_t\omega),z)\\
&-B(z_{W_1}(\theta_t\omega)+z_{W_2}(\theta_t\omega),z_{W_1}(\theta_t\omega)+z_{W_2}(\theta_t\omega))
\\
&-\beta G(z_{W_1}(\theta_t\omega)+z_{W_2}(\theta_t\omega))_x
-r(z_{W_1}(\theta_t\omega)+z_{W_2}(\theta_t\omega))
\end{split}
\end{equation}
with $z(0)=z_0\in H$.
\\
\begin{lemma}\label{lem4.1}
The random evolution equation   (\ref{eq4.0}) has a unique solution for any
initial condition
$z_0$ $\in$ $H$ and  this solution defines a random dynamical system with
respect to the
metric  dynamical system $\theta$  introduced in Section
3 for which
the associated cocycle mapping is defined  by $(t,\omega,z_0)$ $\mapsto$
$z(t,\omega,z_0)$.
\end{lemma}
For the proof of this lemma we can use the fact that equation (\ref{eq4.0})
is quite similar to the Navier--Stokes equation. Although some linear terms
are also present,  similar a priori estimates  can be obtained  to those
in Temam
(\cite{Tem79},  Chapter III)  or Benssousan and Temam \cite{BenTem73}.
Because of
 the properties of the operator $B$ introduced in the previous  section.
Moreover, the random terms appearing inside the coefficients of  equation
(\ref{eq4.0})  are given by stationary processes, so we obtain a random
dynamical
system, see Arnold (\cite{Arn98}, page 58). \\

\begin{remark}
{\rm To see that $t$ $\mapsto$ $ z(t,\omega,z_0)$ is continuous for any
$z_0$ $\in$ $H$
and  $\omega$ $\in$ $\Omega$ we can use   Lemma  III.1.2 in  \cite{Tem79}
since
the solution of  equation  (\ref{eq4.0})   satisfies
\[
\int_0^t\|z(\tau,\omega,z_0)\|^2d\tau<\infty
\]
for any $z_0$ $\in$ $H$.  Indeed, by the chain rule,
\begin{eqnarray}\nonumber
|z(t)|^2 +2\nu\int_0^t\|z(\tau)\|^2d\tau  & \le  & |z_0|^2+2(\beta
c_{G,x}-r)\int_0^t|z(\tau)|^2d\tau
\\
&  & + 2c_B\int_0^t|z(\tau)|\|z(\tau)\||\nabla z_{W_1}(\theta_\tau\omega)  +
\nabla z_{W_2}(\theta_\tau\omega)|d\tau  \nonumber \\
&   &
+ 2\int_0^t\|f(\tau)\|_{-1}\|z(\tau)\|d\tau,  \nonumber
\end{eqnarray}
where $f$ consists of  all the terms in (\ref{eq4.0}) that do not contain $z$.
Then, using
\begin{eqnarray*}
2|z|\|z\||\nabla z_{W_1}+\nabla z_{W_2}|d\tau & \le & \frac{\nu}{2} \|z\|^2
+\frac{2}{\nu} |z|^2 |\nabla z_{W_1}+\nabla z_{W_2}|^2,
\\
\qquad
2\|f\|_{-1} \|z\| & \le &  \frac{\nu}{2} \|z\|^2+\frac{2}{\nu} \|f\|_{-1}^2,
\end{eqnarray*}
the asserted estimate follows by  an application of the Gronwall  inequality.
Moreover, by the properties of the operators $A$ and $B$,  we also have
\[
\int_0^t\|z(\tau,\omega,z_0)\|_{V^\prime}^2d\tau<\infty,\quad z_0\in H.
\]
} \\
\end{remark}

We now define the random isomorphism $i(\omega)$  $:$  $H$ $\to$  $H$
by
$$
i(\omega)a =-(z_{W_1}(\omega)+z_{W_2}(\omega))+a,
$$
for which the inverse isomorphism $i^{-1}(\omega)$ is given by
\[
i^{-1}(\omega)a=(z_{W_1}(\omega)+z_{W_2}(\omega))+a.
\]
Note  that the random variable $i(\omega)a(\omega)$ is tempered for any
tempered
$a(\omega)$.
\begin{lemma}
Let $z(\cdot,\omega,z_0)$ be the solution of (\ref{eq4.0}).  Then the process
\[
u(t,\omega,u_0)=i^{-1}(\theta_t\omega)\circ z(t,\omega,i(\omega)\circ u_0)
\]
solves (\ref{eq3.3}).  In particular, $u$ satisfies the boundary conditions
(\ref{eq3.4}).
\end{lemma}
\begin{proof}   The assertion follows by replacing $z$ by $u$ $-$
$z_{W_1}(\theta_t \omega)$
$-$  $z_{W_2}(\theta_t\omega)$.
\end{proof}

We will now check in the following Lemmata  that the assumptions of the
random fixed
point theorem \ref{theo2.1}  are satisfied.   First, we show that there
exists a tempered
random set ${\cal X}(\omega)$ of (single valued)  random variables that
will be mapped into
itself.
\begin{lemma}\label{l4.4}
 Let  ${\cal X}(\omega)$  be the ball $B(0,\rho(\omega))$ in $H$ with
center zero and
${\cal F}_0$--measurable  radius
\begin{eqnarray*}
\rho(\omega) &=& \bigg(\int_{-\infty}^0 \exp
\bigg((\lambda_1\nu-2\beta C_{G,x}+2r)\tau
\\
& &  \quad + \frac{3c_B^2}{\nu}\int_\tau^0|\nabla z_{W_1}(\theta_s\omega)
+ \nabla z_{W_2}(\theta_s\omega))|^2ds\bigg)\cdot
R(\theta_\tau\omega)d\tau\bigg)^\frac{1}{2},
\end{eqnarray*}
where
\begin{eqnarray*}
R(\omega)&=&3\frac{(C_{G,x}\beta+r)^2}{\nu\lambda_1}|z_{W_1}(\theta_t\omega)+z_{W_2}(\theta_t\omega)|^2\\
&+&
\frac{3c_B^2}{\nu}|z_{W_1}(\theta_t\omega)+z_{W_2}(\theta_t\omega)|^2
|\nabla z_{W_1}(\theta_t\omega)+\nabla z_{W_2}(\theta_t\omega)|^2,
\end{eqnarray*}
and suppose that
\[
\lambda_1\nu+2r>2C_{G,x}\beta +\frac{3c_B^2}{\nu}\EX|\nabla z_{W_1}+\nabla
z_{W_2}|^2,
\]
where $\lambda_1>0$ is the first eigenvalue of the operator $A$.  Then  the
random set ${\cal X}$  is
forward invariant, i.e.,
\[
z(t,\omega,{\cal X}(\omega))\subset {\cal X}(\theta_t\omega),  \qquad t\ge 0.
\]
\end{lemma}
\begin{proof}We have to  estimate $|z|^2$  for which we need  the following
relations
that are a consequence of  Lemma \ref{lem2.3}:
\begin{eqnarray*}
&&
2\langle B(z,z),z\rangle=0,\qquad
2\langle Az,z\rangle=2\nu\|z\|^2\ge \nu \|z\|^2+\nu\lambda_1|z|^2,
\\[2ex]
&&2\beta\langle G(z)_x,z\rangle\le 2\beta C_{G,x}|z|^2,\qquad
2r\langle z,z\rangle =2r|z|^2,
\\[2ex]
&&
2|\langle B(z,z_{W_1}(\theta_t\omega)+z_{W_2}(\theta_t\omega)),z\rangle  |\le
2c_B|z||\nabla z_{W_1}(\theta_t\omega)+\nabla z_{W_2}(\theta_t\omega))|\|z\|
\\
&& \qquad  \qquad  \qquad  \qquad \le \frac{3c_B^2}{\nu}
 |\nabla z_{W_1}(\theta_t\omega)+\nabla
z_{W_2}(\theta_t\omega))|^2|z|^2+\frac{\nu}{3}\|z\|^2,
 \\[2ex]
&&
2\langle B(z_{W_1}(\theta_t\omega)+ z_{W_2}(\theta_t\omega)),z),z\rangle=0,
\\[2ex]
&&
2|\langle B(z_{W_1}(\theta_t\omega)+z_{W_2}(\theta_t\omega),z_{W_1}(\theta_t\omega)
+z_{W_2}(\theta_t\omega)),z\rangle|\\
&&\quad \qquad  \qquad \le
\frac{3c_B^2}{\nu}|z_{W_1}(\theta_t\omega)+z_{W_2}(\theta_t\omega)|^2
|\nabla z_{W_1}(\theta_t\omega)+\nabla z_{W_2}(\theta_t\omega)|^2+\frac{\nu}{3}\|z\|^2,
\\[2ex]
&&
2|\langle\beta G(z_{W_1}(\theta_t\omega)+z_{W_2}(\theta_t\omega))_x-r\langle
z_{W_1}(\theta_t\omega)+z_{W_2}(\theta_t\omega),z\rangle|\\
&&\quad \qquad  \qquad
\le 3\frac{(C_{G,x}\beta+r)^2}{\nu\lambda_1}|z_{W_1}(\theta_t\omega)
+z_{W_2}(\theta_t\omega)|^2+ \frac{\nu\lambda_1}{3}|z|^2.
\end{eqnarray*}
It can be shown by a comparison argument that $|z(t,\omega,z_0)|^2$ is
bounded by a
solution of the affine random differential equation:
\begin{equation}\label{eq4.2}
\begin{split}
\frac{d\zeta}{dt}&+(\lambda_1\nu-2\beta C_{{G,x}}
+2r-\frac{3c_B^2}{\nu}(|\nabla z_{W_1}(\theta_t\omega)+
\nabla z_{W_2}(\theta_t\omega))|^2))\zeta =R(\theta_t\omega) \\[2ex]
&  \qquad \qquad \qquad\zeta(0)=|z_0|^2,
\end{split}
\end{equation}
for which  the solution is given by a variation of constant formula.
A  direct calculation verifies   that $t$ $\to$ $\rho^2(\theta_t\omega)$
is a solution of equation (\ref{eq4.2}) with initial value
$\zeta(0)$ $=$ $\rho^2(\omega)$, which  means that $\rho^2$ is a random
fixed point of
(\ref{eq4.2}). It thus  follows that   $z(t,\omega,z_0)$ $\in$ ${\cal
X}(\theta_t\omega)$
whenever  $z_0$ $\in$ ${\cal X}(\omega)$.
\end{proof}
We note that  the random variable $\rho$ is tempered (see  \cite{Schm97a} ,
page 110),
so  any selector contained in ${\cal X}$ is also tempered.\\

It remains to check that the contraction condition (\ref{eq2.2}) of the
random fixed point theorem
holds.
\begin{lemma}\label{l4.5}
Suppose that
\begin{equation}\label{eq4.4}
\begin{split}
&-\nu\lambda_1
+2\beta C_{G,x}-2r+\frac{3c_B^2}{\nu}\EX|\nabla z_{W_1}+\nabla z_{W_2}|^2+\\
&\frac{2c_B^2}{\nu^2}(1+2\lambda_1^\frac{1}{2}C_{G,x}\beta
)\EX\rho^2
+\frac{c_B^2}{\nu}\EX\rho^4+
\frac{c_B^2}{\nu}\EX|\nabla z_{W_1}+\nabla z_{W_2}|^4 +\frac{2\EX R}{\nu}<0,
\end{split}
\end{equation}
where $R$ was defined in Lemma \ref{l4.4}. Then the contraction condition
(\ref{eq2.2})
is fulfilled.
\end{lemma}
\begin{proof}
It follows immediately  from Lemma \ref{lem2.3} with $z_1$, $\,z_2$ $\in$
$V$
that
\begin{eqnarray}
|\langle B(z_1,z_1)-B(z_2,z_2),z_1-z_2\rangle|
& = &  |\langle B(z_1,z_1),z_2\rangle+\langle B(z_2,z_1),z_2\rangle|
\nonumber \\
& \le &
|\langle B(z_1-z_2,z_1),z_1-z_2\rangle| \nonumber\\
& \le &
c_B|z_1-z_2|\|z_1\|\|z_1-z_2\|. \nonumber
\end{eqnarray}
Set $z_1=z(t,\omega,h_1)$, $z_2=z(t,\omega,h_2)\in H$ and $\delta z$ $=$ $z_1-z_2$.
By the chain rule we obtain
\begin{eqnarray*}
|\delta z(1)|^2+2\nu\int_0^1\|\delta z\|^2d\tau & \le &|h_1-h_2|^2+
\int_0^1(2C_{G,x}\beta|\delta z|^2-2r|\delta z|^2\\
&  & +2c_B|\nabla z_{W_1}(\theta_s\omega)+\nabla z_{W_2}(\theta_s\omega))|\|\delta
z\||\delta z|
\\
&  & +
2c_B\|\delta z\|\|z_1\||\delta z|)ds
\end{eqnarray*}
for $h_1$, $\,h_2$ $\in$ $ {\cal X}(\omega)$.  It  then follows from
\begin{eqnarray*}
&&
2c_B|\nabla z_{W_1}(\theta_s\omega)+\nabla z_{W_2}(\theta_s\omega))|\|\delta
z\||\delta z|
 \le  \frac{3c_B^2}{\nu}|\nabla z_{W_1}+\nabla z_{W_2}|^2|\delta
z|^2+\frac{\nu}{2}\|\delta z\|^2
\\[2ex]
&&
2c_B\|\delta z\|\|z_1\||\delta z|  \le \frac{\nu}{2}\|\delta
z\|^2+\frac{2c_B^2}{\nu}\|z_1
\|^2|\delta z|^2,\qquad\nu\|z\|^2\ge \nu\lambda:1|z|^2
\end{eqnarray*}
that
\[
\frac{1}{2}\EX\sup_{h_1\not=h_2\in {\cal X}(\omega)}
\log\frac{|z(1,\omega,h_1)-z(1,\omega,h_2)|^2}{|h_1-h_2|^2}
\]
is less than or equal to the expression on the left hand side of inequality
(\ref{eq4.4}).
Here we have used the fact that
\begin{eqnarray*}
&&
\frac{2c_B^2}{\nu}\int_0^1\|z_1(t,\omega,z_0)\|^2d\tau  \le
\frac{2c_B^2}{\nu^2}
\bigg(\rho(\omega)^2 +(2C_{{G,x}}\beta-2r)\int_0^1\rho(\theta_\tau\omega)^2d\tau
\\
&& \qquad +  \frac{c_B^2}{\nu}\int_0^1\rho(\theta_\tau\omega)^4d\tau
\\
&&\qquad  +\frac{c_B^2}{\nu}\int_0^1|\nabla z_{W_1}(\theta_\tau\omega)+\nabla
z_{W_2}(\theta_\tau\omega)|^4d\tau
+\frac{2\int_0^1R(\theta_\tau\omega)d\tau}{\nu}\bigg)
\end{eqnarray*}
for $z_0$ $\in$ ${\cal X}(\omega)$.
\end{proof}
\begin{remark}
{\rm
Note that the assumption of Lemma \ref{l4.5} is sufficient for (\ref{l4.4}).
}
\end{remark}
The crucial point for the assumptions of the last lemma is to show  at
least for large $\nu$
and for small
${\rm tr}_{U}Q^{W_1}$ and ${\rm tr}_{{W_2}} Q^V$ that the
random variables
$\rho^2$ and $\rho^4$ have finite and sufficiently small expectations.
  In addition, the finiteness of the expectation of these
random variables ensures that (\ref{eq2.2.1}) is satisfied.

\begin{lemma}
 The expectations of $\rho^2$ and $\rho^4$ are sufficiently small  when
$\nu$  is
 sufficiently large and  ${\rm tr}_{U}Q^{W_1}$ and ${\rm
tr}_{{V}} Q^{W_2}$ are
sufficiently small.
\end{lemma}
We   give only a brief comment on the proof of this very technical lemma.
The essential ingredient is that $z_{W_2}$ and $\nabla z_{W_1}$ are
Gaussian
processes, so
$R$ has finite moments of arbitrary order and
\begin{equation}\label{eq4.6}
E\exp\left(\alpha_1\int_0^t|\nabla z_{W_1}|^2d\tau\right)\le e^{\alpha_2t},
\end{equation}
and similarly for $\nabla z_{W_2}$. The constant $\alpha_1>0$ depends of the
data  $\nu$, $\lambda_1$, $\cdots$ of the problem. The assertion
of the Lemma follows if  $\alpha_2$ is sufficiently small,  which can be
controlled
by the traces of $Q^{W_2}$ and $ Q^{W_1}$.
Finally,  to obtain  a finite dimensional version of the estimate
(\ref{eq4.6}) we refer
to Hasminski{\v{i}} (\cite{Has80}, page 37,  Lemma 7.2),  where we need the
main
assumption
\[
E\,{\rm tr_{H}}(\nabla z_{W_1}(\theta_{t_1}\omega)\otimes \nabla
z_{W_1}(\theta_{t_2}\omega))
\le c\,e^{-\lambda_1|t_1-t_2|}
\]
for an appropriate constant $c$. The variable $z_{W_2}$ can be handled similarly.\\

Summarising,  we have
\begin{theorem}
Suppose that the assumption of Lemma \ref{l4.5} is satisfied
and let   $z^\ast$ be the random fixed point of the(transformed)
random dynamical system generated by (\ref{eq4.0}).
Then there exists a random fixed point for (\ref{eq3.3}) that  attracts
the states of the phase space exponentially fast.
\end{theorem}

Indeed, the random variable $u$ that  generates  an exponentially stable
stationary solution is given by
\[
u^\ast(\omega)=z^\ast(\omega)+z_{W_2}(\omega)+z_{W_1}(\omega).
\]

\section{Discussions}

We have shown that, under suitable conditions on the random
forcing, random boundary conditions, viscosity, Ekman constant
and Coriolis parameter,   all quasigeostrophic motion
approach a unique stationary
state exponentially fast as time goes to infinity.
In   deterministic systems a high level of stability is obtained
when there is an exponential attractor which attracts trajectories
exponentially fast. In some situations this attractor is a single point
(point attractor)
which describes the laminar behavior of the flow.
We are looking for such   stability in the case of quasigeostrophic fluid motion
under random perturbations.
In particular, we find a random attractor which is defined by a single
random variable. This  random variable attracts  all other quasigeostrophic
motion exponentially fast. This random variable corresponds to a unique
invariant measure, which is the Dirac measure with the random variable
as the random mass point; see \cite{Arn98}.
The corresponding stationary Markov measure is
the expectation of this random Dirac measure.

\end{document}